\documentclass[runningheads]{llncs}
\usepackage{graphicx}
\usepackage{amssymb}
\usepackage{url}
\usepackage{algorithm}
\usepackage{algpseudocode}
\usepackage{fontawesome5}

\usepackage{amsmath}

\newcommand{\mK}{\mathsf{K}}

\usepackage{multirow}

\usepackage{color}

\begin{document}
\title{A Bayesian Approach for Simultaneously Radial
Kernel Parameter Tuning in the Partition of
Unity Method}
\titlerunning{A Bayesian Approach for Shape Parameter Tuning in RBF -PUM}
%
\author{Roberto Cavoretto\orcidID{0000-0001-6076-4115} \and
 Alessandra De Rossi\orcidID{0000-0003-1285-3820} \and
 Sandro Lancellotti\orcidID{0000-0003-4253-3561} \and
 Federico Romaniello\orcidID{0000-0003-1166-3179}\faIcon{envelope}}
\authorrunning{R. Cavoretto et al.}
%
\institute{Department of Mathematics \lq\lq Giuseppe Peano\rq\rq, University of Torino, via Carlo Alberto 10, 10123 Torino, Italy \\
\email{roberto.cavoretto@unito.it, alessandra.derossi@unito.it, sandro.lancellotti@unito.it, \faIcon{envelope}federico.romaniello@unito.it}}
\maketitle              
\begin{abstract}
In this paper, Bayesian optimisation is used to simultaneously search the optimal values of the shape parameter and the radius in radial basis function partition of unity interpolation problem. It is a probabilistic iterative approach that models the error function with a step-by-step self-updated Gaussian process, whereas partition of unity leverages a mesh-free method that allows us to reduce cost-intensive computations when the number of scattered data is very large, as the entire domain is decomposed into several smaller subdomains of variable radius. Numerical experiments on the scattered data interpolation problem show that the combination of these two tools sharply reduces the search time with respect to other techniques such as the leave one out cross validation. 
\end{abstract}

\section{Introduction}
Meshfree methods are popular tools for solving interpolation problems and one of the most used is the Partition of Unity Method (PUM). Here, it is implemented using Radial Basis Functions (RBFs) or radial kernels as local approximants, since PUM is a very efficient tool in scattered data interpolation \cite{Sandro,Wendland05}. The main disadvantage of radial kernel-based method is the computational cost associated with the solution of (usually) large linear systems; PUM are instead used to efficiently split the data into smaller subdomains (or balls) of radius $\delta$. Such a method is obtained by a weighted sum of some RBFs depending on a {\sl shape parameter} $\varepsilon$ in each subdomain \cite{cav21a}.

The aim of this work is to apply a well-known statistical technique, called {\sl Bayesian Optimization} (BO) \cite{Practical}, to simultaneously search the optimal $(\varepsilon,\delta)$ values for each subdomain in RBF-PUM. This technique, developed in machine learning  for optimisation of black-box or difficult-to-evaluate functions, can be used in hyperparameter tuning problems to avoid computation and evaluation of the interpolant for those parameters that are far from being optimal, leading to a significant reduction of the computational time. To complete this analysis, we compare BO results with those obtained by using the Leave One Out Cross Validation (LOOCV) scheme \cite{Rippa}.

The paper is organised as follows. In Section \ref{rbfpum}, RBF-PUM interpolation and LOOCV method are briefly stated. In Section \ref{bo}, BO Gaussian processes and acquisition functions are described. In Section \ref{exp}, numerical experiments are analysed and discussed. Section \ref{conclusion} concludes the paper.

\section{RBF-PUM Interpolation}\label{rbfpum} 

In this section, we introduce the interpolation problem and the basic theory on RBF-PUM, highlighting the reasons that inspired this contribution. 

\subsection{The RBF Method}
Given a set of distinct data $ X = \{  \boldsymbol{x}_i, i = 1,  \ldots , n \}$ arbitrarily distributed on a domain $ \Omega \subseteq \mathbb{R}^{d}$ with an associated set $ F = \{ f_i = f(\boldsymbol{x}_i) ,i=1, \ldots, n \}$ of data values obtained by sampling some, possibly unknown, function  $f: \Omega \rightarrow \mathbb{R}$ at the nodes $ \boldsymbol{x}_i$, the \emph{scattered data interpolation problem} consists in finding a function $P_f: \Omega \rightarrow \mathbb{R}$ such that it matches the measurements at the corresponding locations, i.e. $P_f\left( \boldsymbol{x}_i\right)=f_i$, $i=1, \ldots, n$.

We now suppose to have a univariate function $ \varphi: [0, \infty) \to \mathbb{R}$, known as RBF, which depends on a shape parameter $\varepsilon > 0$ providing, for $\boldsymbol{x},\boldsymbol{z}
\in \Omega$, the real symmetric strictly positive definite kernel 
\begin{equation*} 
\kappa_\varepsilon(\boldsymbol{x},\boldsymbol{z}) = \varphi(\varepsilon ||\boldsymbol{x}-\boldsymbol{z}||_2): = \varphi(\varepsilon r).
\end{equation*}

The kernel-based interpolant $P_f$ can be written as
\begin{equation*}
P_f\left( \boldsymbol{x}\right)= \sum_{k=1}^{n} c_k \kappa_\varepsilon \left( \boldsymbol{x} , \boldsymbol{x}_k  \right), \quad \boldsymbol{x} \in \Omega,
\end{equation*}
whose coefficients are the solution of the linear system
\begin{equation} \label{linsys}
\mK \boldsymbol{c} = \boldsymbol{f},
\end{equation}
where $  \boldsymbol{c}= \left(c_1, \ldots,
c_n\right)^{\intercal}$, $  \boldsymbol{f} =\left(f_1, \ldots , f_n\right)^{\intercal}$, and $\mK_{ik}= \kappa_\varepsilon \left( \boldsymbol{x}_i , \boldsymbol{x}_k  \right)$, $i,k=1, \ldots, n$. Since $\kappa_\varepsilon$ is a symmetric and strictly positive definite kernel, the system \eqref{linsys} has exactly one solution \cite{Fasshauer15}. 

\subsection{The PUM Scheme and the LOOCV Technique}
It is well-known that inverting the kernel interpolation matrix in \eqref{linsys} might be computationally expensive when the amount of data highly increases. An effective method to overcome this problem is to partition the open and bounded domain $\Omega$ into $m$ overlapping subdomains, i.e $\Omega \subseteq \bigcup_{j=1}^m\Omega_j$, and hence to split the interpolation problem locally in each one of them.

The PU covering consists of overlapping balls of radius $\delta$ whose centres are the grid data $P=\{\tilde{\boldsymbol{x}}_k$, $k= 1,\ldots,m\}$.

In \cite{Fasshauer} it is shown that when the nodes are nearly uniformed distributed, $m$ is a suitable number of PU subdomains on $\Omega$ if ${n}/{m} \approx 2^d$. Then, the covering property is satisfied by taking the radius $\delta$ such that
\begin{equation*} 
      	\delta \geq \frac{\displaystyle 1 }{\displaystyle m^{1/d}}.
      	\label{PU_radius}
\end{equation*}

The PUM solves a local interpolation problem on each subdomain and constructs the global approximant by gluing together the local contributions using weights. To achieve that, we need those weights are a family of compactly supported, non-negative, continuous functions $w_j$, with $\text{supp}\left(w_j\right) \subseteq \Omega_j$, such that
\begin{equation*}
\sum_{j=1}^{m} w_j\left(\boldsymbol{x}\right) = 1, \quad \boldsymbol{x} \in \Omega.
\end{equation*}

Once we choose the partition of unity $ \{ w_j \}_{j=1}^{m}$, the global interpolant is formed by the weighted sum of $m$ local approximants $P_f^j$, i.e.
\begin{equation*}
P_f\left( \boldsymbol{x}\right)= \sum_{j=1}^{m} P_f^j\left( \boldsymbol{x} \right) w_j \left( \boldsymbol{x}\right) = \sum_{j=1}^{m} \left( \sum_{k=1}^{n_j} c_k^j \kappa_{\varepsilon,\delta}(\boldsymbol{x} ,  \boldsymbol{x}^j_k) \right)  w_j \left( \boldsymbol{x}\right), \quad \boldsymbol{x} \in \Omega,
\label{intg}
\end{equation*}
where $n_j=|\Omega_j|$ and $\boldsymbol{x}^j_k \in X_j = X \cap \Omega_j$, with $k=1, \ldots, n_j$.

We remark that the accuracy of the fit strongly depends on the choices of the shape parameter and the radius, see e.g. \cite{cav21,cav22,Fornberg-Wright04,Larsson-Fornberg05,lin22}.

Rippa in \cite{Rippa} proposed the LOOCV for the search of the optimal value of the RBF shape parameter $\varepsilon$. This technique can be extended to simultaneously search the optimal values for $\delta$ and $\varepsilon$ and applied in each PU subdomain $\Omega_j$. It consists in evaluating in every subdomain $\Omega_j$, for each $(\varepsilon,\delta)$ and $ k \in \lbrace 1, \dots, n_j \rbrace$, the interpolation error $e_{j,k}(\varepsilon,\delta)$ at the point $x_k$ of the RBF interpolant $P_f^{j,k}$  fitted on the data set $X_k = X \setminus \lbrace x_k \rbrace$ and the data values  $F_k = F \setminus \lbrace f_k \rbrace$. To avoid this computation, the error can be computed by the rule:
\begin{equation*}
    e_{j,k}(\varepsilon,\delta) = \frac{c_k^j}{(\mK_{kk}^{-1})_j}.
\end{equation*}

Hence, the optimal value $(\varepsilon,\delta)_j^*$ is the one that minimises the error function $Er(\varepsilon,\delta)_j$ defined as follows:
\begin{equation}\label{error_function}
Er(\varepsilon,\delta)_j = \max_{k=1, \dots, n_j} \Bigg|   \frac{c_k^j}{(\mK_{kk}^{-1})_j} \Bigg|.
\end{equation}

The LOOCV technique is then formalised by imposing to evaluate the error function also for that value in the discrete set that does not lead to a good result. A drawback of this scheme is the impossibility, in general, to attain the global minimum due to the discrete research. In this work we propose a faster procedure that conducts a continuous research in the parameters space in order to obtain a better approximation of $(\varepsilon,\delta)_j^*$.

\section{Bayesian Optimisation}
\label{bo}
When it comes to find a global maximiser for an unknown or difficult-to-evaluate function $f$ on some bounded set $X$, the Bayesian Optimisation (BO) \cite{Mockus_1978} is an approach to carry out the research. Very popular in machine learning, BO is an iterative technique that is based on exploiting all the available resources. It consists in building a probabilistic model of $f$, called  {\sl surrogate model}, and using it to help directly the sampling point in $X$, by means of an acquisition function, where the target function will be evaluated. Once an iteration is made, the distribution is updated and then used in the next iteration. Though there is a computation for the selection of the next point to evaluate, when evaluations of $f$ are expensive, the computation of a better point is motivated by reaching the maximum in a few iterations, as in the case of the error function of some expensive training machine learning algorithms like multi-layer neural networks. Hereinafter, we briefly review the BO technique \cite{Brochu}.

 A Gaussian Process (GP) is a collection of random variables such that any subsets of these have a joint Gaussian distribution. Then GPs are completely specified by a mean function $m : \mathcal{X} \rightarrow \mathbb{R}$ and a positive definite covariance function $k : \mathcal{X} \times \mathcal{X} \rightarrow \mathbb{R}$ (see \cite{Rasmussen}). Moreover, they are the most common choice for the surrogate model for BO due to the low evaluation cost and to the ability of incorporating prior beliefs about the objective function. When modeling the target function with a GP as $f(\textbf{x}) \sim \mathcal{GP} \big( m(\textbf{x}), k(\textbf{x}, \textbf{x}') \big)$, we impose that
\begin{align*}
    \mathbb{E}\big [ f(\textbf{x}) \big] = m(\textbf{x}), \qquad
    \mathbb{E}\big [ \big( f(\textbf{x}) - m(\textbf{x}) \big)  \big( f(\textbf{x}') - m(\textbf{x}') \big) \big] = k(\textbf{x}, \textbf{x}').
\end{align*}
In the matter of making a prediction given by some observations, the assumption of joint Gaussianity allows retrieving the prediction using the standard formula for mean and variance of a conditional normal distribution. Hence, suppose to have $s$ observation $\textbf{f} = (f(\textbf{x}_1), \dots, f(\textbf{x}_s))^{\intercal}$ on $\textbf{X} = (\textbf{x}_1, \dots, \textbf{x}_s)^{\intercal}$ and a new point $\bar{\textbf{x}}$ on which we are interested in having a prediction of $\bar{f} = f(\bar{\textbf{x}})$. The previous observations $\textbf{f}$ and the predicted value $f(\bar{\textbf{x}})$ are jointly normally distributed:
$$
Pr\begin{pmatrix} \begin{bmatrix} \textbf{f} \\ f(\bar{\textbf{x}}) \end{bmatrix} \end{pmatrix} = \mathcal{N} \begin{bmatrix}
\begin{bmatrix} \mu(\textbf{X}) \\ \mu(\bar{\textbf{x}}) \end{bmatrix}, 
\begin{bmatrix} 
K(\textbf{X},\textbf{X}) \ \ K(\textbf{X}, \bar{\textbf{x}}) \\
K(\textbf{X}, \bar{\textbf{x}})^{\intercal} \ \ k(\bar{\textbf{x}}, \bar{\textbf{x}}) \\
\end{bmatrix}\end{bmatrix},
$$
where $K(\textbf{X},\textbf{X})$ is the $s \times s$ matrix with $(i, j)$-element $k(\textbf{x}_i, \textbf{x}_j)$, and
$K(\textbf{X}, \bar{\textbf{x}})$ is a $s \times 1$ vector whose element $i$ is given by $k(\textbf{x}_i, \bar{\textbf{x}})$. Since $Pr(f(\bar{\textbf{x}}) |  \textbf{f})$ must also be normal, it is also possible to estimate the distribution, the mean and the covariance, for any point in the domain. When data points and data values retrieved by the evaluation of the target function are fed to the model, they induce a posterior distribution over functions which is used for the next iteration as a prior. We point out that if a function is modelled by a GP, when we observe a value, we are observing the random variable associated to the point.



\emph{An acquisition function} $a: \mathcal{X} \rightarrow \mathbb{R}$ is a function used to determine the next evaluation point for the objective function. This chosen point maximises the acquisition function, and its evaluation by the objective function is used to update the surrogate model. An acquisition function is defined such that high acquisition corresponds to potentially high values of the objective function. There exists a trade-off between exploration and exploitation in the selection of an acquisition function: exploration means selecting points where the uncertainty is high, that is, far from the already evaluated points; exploitation, on the contrary, means selecting those points close to those already evaluated by the objective function. The most common acquisition functions are:

\begin{itemize}
    \item \textbf{Probability of Improvement}, which maximises the probability of improvement over the best current value;
    \item \textbf{Expected Improvement}, which maximises the expected improvement over the current best;
    \item \textbf{GP Upper Confidence Bound}, which minimises the cumulative regret.
\end{itemize}
The acquisition function used in this work is the \lq\lq Expected Improvement\rq\rq\ \cite{EI} as it considers not only the probability of improvement of the candidate point w.r.t. the previous maximum, but also the magnitude of this improvement.

Suppose that after a number of iterations the current maximum of the objective function is $f(\hat{\textbf{x}})$. Given a new point $\textbf{x}$, the Expected Improvement acquisition function computes the expectation of improvement $f(\textbf{x}) - f(\hat{\textbf{x}})$  over the part of the normal distribution that is above the current maximum:
\begin{equation}
\label{EI_eq}
    EI(\textbf{x}) = \int_{f(\hat{\textbf{x}})}^\infty \big( f^*(\textbf{x}) - f(\hat{\textbf{x}}) \big) 
    \frac{1}{\sqrt{2 \pi} \sigma(\textbf{x})} e^{-\frac{1}{2}  [(f^*(\textbf{x}) - \mu(\textbf{x}))/\sigma(\textbf{x})]^2} df^*(\textbf{x}),
\end{equation}
where $f^*(\textbf{x})$, $\mu(\textbf{x})$ and $\sigma(\textbf{x})$ represent the predicted value by the surrogate model, the expected value and the variance of $\textbf{x}$, respectively.
Solving integral \eqref{EI_eq} leads to the following closed form for the evaluation of the Expected Improvement:
\begin{equation*} 
    EI(x) = 
    \begin{cases}  
            (\mu(\textbf{x}) - f(\hat{\textbf{x}})) \Phi(Z) + \sigma(\textbf{x)} \phi(Z), & \text{ if } \sigma(\textbf{x}) >0, \\
            0, & \text{ if } \sigma(\textbf{x}) =0, 
    \end{cases}
\end{equation*}
where $Z =\frac{\mu(\textbf{x})-f(\hat{\textbf{x}})}{\sigma(\textbf{x})}$, while $\phi$ and $\Phi$ are the Probability Density Function and Cumulative Distribution Function of the standard normal distribution $\mathcal{N}(0, 1)$. 

\section{Numerical Experiments} \label{exp}
 A sketch of the proposed method is shortly described in Algorithm \ref{alg33}.
 \begin{algorithm}
\caption{BO-PUM}
\hspace*{\algorithmicindent} \textbf{Input:} Data points $X$, data values $F$.
\begin{algorithmic}
\State Construct a partition of unity with centers generated as an equispaced grid;
\State In each subdomain, find the value for the radius that ensures the minimum density;
\State Apply BO in each subdomain to find $(\varepsilon,\delta)$;
\State For each subdomain construct a local RBF approximant;
\State Create a global RBF-PUM approximant from the local ones.
      \end{algorithmic}
\hspace*{\algorithmicindent} \textbf{Output:} RBF-PUM approximation.
\label{alg33}
\end{algorithm}
Numerical experiments are presented below to find the best $(\varepsilon,\delta)$ that minimise the error function for each subdomain.
We compare the results obtained by applying BO and LOOCV in the RBF-PUM. To apply the BO for the search of the optimal parameters $(\varepsilon,\delta)$, we suppose that the objective function to maximise is the Maximum Absolute Error (MAE) of the RBF interpolant, changed of sign because the BO is a maximisation process as described in Section \ref{bo}. We use a variable number of points for the training set and a fixed one of $1000$ points for the test set. During the BO, for each subdomain, after determining the points belonging to it,  we further divide them in a sub-training and sub-validation set to allow the evaluation of training error. After identifying the best parameter pairs for each subdomain and each optimiser, a PUM interpolant is trained on the training set for each optimiser with the parameters found.
For each subdomain, the search space is $(0, 20] \times [\delta_{min}, 2 \delta_{min}]$, where $\delta_{min}$ is the radius value that ensures a minimum density in the subdomain. The LOOCV tries each possible combination between $500$ equally spaced points in  $(0, 20]$ and $30$ equally spaced points in $[\delta_{min}, 2 \delta_{min}]$. On the other hand, BO performs $5$ random steps plus at most $25$ Bayesian steps in the search space. The iterative process stops if  the desired tolerance $\tau$ is reached.
We perform the experiments on $3$ different sizes of random data in the interval $[0, 1]^2$ using the RBFs 
\begin{align*}
&\varphi_1(\varepsilon r) = e^{-\varepsilon^2 r^2} & (\mbox{Gaussian $C^{\infty}$}), \\ 
&\varphi_2(\varepsilon r)= e^{-\varepsilon r}(1 + 3 \varepsilon r + \varepsilon^2 r^2) & (\mbox{Mat$\acute{\text{e}}$rn $C^4$}),
\end{align*}
and the test function 
\begin{align*}
    f(x_1,x_2) =  2 \cos(10 x_1) \sin(10 x_2) + \sin(10 x_1 x_2).
\end{align*}

Results are shown in Table \ref{tab:f}. As it can be seen in this table, the errors obtained using the different techniques are similar in all cases, while the runtime of BO is always lower than LOOCV. It should also be noted that, for the Gaussian kernel, as the number of points increases, the execution time of the BO decreases. This fact is due to the high density of the space when a greater number of points is considered.
 In particular, when this happens, there are denser subdomains and thus better accuracy and fewer BO iterations are needed to satisfy the tolerance $\tau$.

\begin{table}[]
    \centering
    
\begin{tabular}{|*{7}{p{17mm}|}}

\cline{4-7}
    \multicolumn{3}{c}{} & \multicolumn{2}{|c}{Gaussian kernel} &\multicolumn{2}{|c|}{Matérn kernel}\\
    \hline   
$N$ &   optimizer  &  $\tau$ &  time (s)  & MAE  & time (s) & MAE \\
\hline
\multirow{3}{*}{8000} & LOOCV &        & 3.95e+03 & 1.87e-04 & 4.06e+03 &  1.68e-03 \\
\cline{2-7}
\multirow{3}{*}{} &    \multirow{2}{*}{BO} & 1e-04 & 1.18e+01 & 7.83e-05 & 3.36e+02 &  4.56e-04 \\
\cline{3-7}
\multirow{3}{*}{} &   \multirow{2}{*}{}  & 1e-05 &  4.19e+01 & 8.84e-06 & 8.37e+02 &  4.96e-04 \\
\hline
\multirow{3}{*}{4000} & LOOCV &        & 1.96e+03 & 1.34e-05 & 2.01e+03 &  6.26e-03 \\
\cline{2-7}
\multirow{3}{*}{} &    \multirow{2}{*}{BO} & 1e-04 & 1.70e+01 & 6.91e-05 & 4.19e+02 &  1.08e-03 \\
\cline{3-7}
\multirow{3}{*}{} &    \multirow{2}{*}{}  & 1e-05 & 1.07e+02 & 4.88e-05 & 4.69e+02 &  3.04e-03 \\
\hline
\multirow{3}{*}{2000} & LOOCV &      & 9.73e+02 & 9.07e-03 & 1.01e+03 &  1.87e-02 \\
\cline{2-7}
\multirow{3}{*}{} &   \multirow{2}{*}{BO} & 1e-04 & 4.01e+01 & 8.67e-05 & 2.81e+02 &  1.24e-02 \\
\cline{3-7}
\multirow{3}{*}{} &    \multirow{2}{*}{}  & 1e-05 & 2.04e+02 & 1.61e-04 & 2.71e+02 &  8.89e-03 \\
\hline
\end{tabular}
    \caption{Computational time and MAE using LOOCV and BO optimiser for Gaussian and Matérn kernels and different number $N$ of random points in $[0,1]^2$. Two tolerances $\tau$ for the BO training error are used.}
    \label{tab:f}
\end{table}

\section{Conclusions}\label{conclusion}
As can be seen from the numerical experiments in Section \ref{exp}, when comparing the parameters search in the case of interpolation with LOOCV and BO, it becomes evident that the error values generally hover around the same magnitude. However, the key distinguishing factor between the two methods is the computational time, wherein the BO approach significantly reduces the time required, often by an order of magnitude. 

\subsubsection*{Acknowledgments}  

The authors sincerely thank the reviewers for the careful reading and valuable comments on the paper. This research has been accomplished within the RITA \lq\lq Research ITalian network on Appro\-xi\-ma\-tion\rq\rq\ and the UMI Group TAA \lq\lq Approxi\-ma\-tion Theory and Applications\rq\rq. This work has been supported by the INdAM--GNCS 2022 Project \lq\lq Computational methods for kernel-based approximation and its applications\rq\rq, code CUP$\_$E55F22000270001, and by the Spoke 1 ``FutureHPC \& BigData'' of the ICSC--National Research Center in "High-Perfor\-man\-ce Computing, Big Data and Quantum Computing", funded by European Union -- NextGenerationEU. Moreover, the work has been supported by the Fondazione CRT, project 2022 \lq\lq Modelli matematici e algoritmi predittivi di intelligenza artificiale per la mobilit$\grave{\text{a}}$ sostenibile\rq\rq.


%
%
%

\begin{thebibliography}{00}



\bibitem{Brochu} E. Brochu, V.M. Cora, N. De Freitas, A tutorial on Bayesian optimization of expensive cost functions, with application to active user modeling and hierarchical reinforcement learning, 2010, arXiv:1012.2599 

\bibitem{cav21a} R. Cavoretto, Adaptive radial basis function partition of unity interpolation: A bivariate algorithm for unstructured data, J. Sci. Comput. 87  (2021) 41.

\bibitem{Sandro} R. Cavoretto, A. De Rossi, S. Lancellotti, E. Perracchione, Software implementation of the partition of unity method, Dolomites Res. Notes Approx. 15 (2022) 35--46. 

\bibitem{cav21} R. Cavoretto, A. De Rossi, M.S. Mukhametzhanov, Ya.D. Sergeyev, On the search of the shape parameter in radial basis functions using univariate global optimization methods, J. Global Optim. 79 (2021) 305--327.

\bibitem{cav22} R. Cavoretto, A. De Rossi, A. Sommariva, M. Vianello, RBFCUB: A numerical package for near-optimal meshless cubature on general polygons, Appl. Math. Lett. 125 (2022) 107704.

\bibitem{Fasshauer} G.E. Fasshauer, Meshfree Approximation Methods with MATLAB, World Scientific, Singapore, 2007.

\bibitem{Fasshauer15}
G.E. Fasshauer, M.J. McCourt, Kernel-based Approximation Methods Using MATLAB, World Scientific, Singapore, 2015.

\bibitem{Fornberg-Wright04} 
B. Fornberg, G. Wright, Stable computation of multiquadrics interpolants for all values of the shape parameter, Comput. Math. Appl. 47 (2004) 497--523.
 


\bibitem{EI} D.R. Jones, M. Schonlau, W.J. Welch, Efficient Global Optimization of Expensive Black-Box Functions, J. Global Optim. 13 (1998) 455--492.


\bibitem{Larsson-Fornberg05} 
E. Larsson, B. Fornberg, Theoretical and computational aspects of multivariate interpolation with increasingly flat radial basis functions, Comput. Math. Appl. 49 (2005) 103--130.

\bibitem{lin22} L. Ling, F. Marchetti, A stochastic extended Rippa’s algorithm for LpOCV, Appl. Math. Letters 129 (2022) 107955.



\bibitem{Mockus_1978} J. Mockus, V. Tiesis, A. Zilinskas, The application of Bayesian methods for seeking the extremum, Towards Global Optimization 2 (1978) 117--129.


\bibitem{Rasmussen} C.E. Rasmussen, C. Williams, Gaussian Processes for Machine Learning, MIT Press, 2006.

\bibitem{Rippa} S. Rippa, An algorithm for selecting a good value for the parameter $c$ in radial basis function interpolation, Adv. Comput. Math. 11 (1999) 193--210.

\bibitem{Practical} J. Snoek, H. Larochelle, R.P. Adams, Practical Bayesian Optimization of Machine Learning Algorithms, Advances in Neural Information Processing Systems 25 (2012) 2960--2968.


\bibitem{Wendland05}
H. Wendland, Scattered Data Approximation, Cambridge Monogr. Appl. Comput. Math., vol. 17, Cambridge Univ. Press, Cambridge, 2005.


\end{thebibliography}
%

\end{document}